\documentclass[12pt]{amsart}

\newtheorem{theorem}{Theorem}[section]

\newtheorem{corollary}[theorem]{Corollary}

\newtheorem{lemma}[theorem]{Lemma}

\theoremstyle{definition}

\newtheorem{definition}[theorem]{Definition}

\newtheorem{remark}[theorem]{Remark}

\newtheorem{example}[theorem]{Example}

\theoremstyle{parrafo}

\numberwithin{equation}{section}
\begin{document}

\title[Derivative of the maximal function]{Functions of bounded variation, the derivative
 of the one dimensional maximal function,
and applications to inequalities}

\author[J. M. Aldaz, J. P\'erez L\'azaro]{J. M. Aldaz and J. P\'erez L\'azaro}
\address{Departamento de Matem\'aticas y Computaci\'on,
Universidad  de La Rioja, 26004 Logro\~no, La Rioja, Spain.}
\email{aldaz@dmc.unirioja.es}
\address{Departamento de Matem\'aticas e Inform\'atica,
Universidad P\'ublica de Navarra, 31006 Pamplona, Navarra, Spain.}
\email{francisco.perez@unavarra.es}


%\thanks{2000 {\em Mathematical Subject Classification.}
%42B25, 26A84}

\thanks{Both authors were partially supported by Grant BFM2003-06335-C03-03 of the
D.G.I. of Spain}

\thanks{The second named author thanks the University of La Rioja for its
hospitality.}

%\subjclass{}


%\date{}

%\dedicatory{}

%\commby{}

%%% ----------------------------------------------------------------------
\begin{abstract} We prove that if $f:I\subset \Bbb R\to \Bbb R$ is of bounded variation, then the uncentered
maximal function $Mf$ is absolutely continuous, and its derivative
satisfies the sharp inequality $\|DMf\|_{L^1(I)}\le |Df|(I)$. This
allows us obtain, under less regularity, versions of classical
inequalities involving derivatives.
\end{abstract}

%%% ----------------------------------------------------------------------

\maketitle

%%% ----------------------------------------------------------------------

\section
{Introduction.}

\markboth{J. M. Aldaz, J. P\'erez L\'azaro} {Derivative of the
maximal function}

The study of the regularity properties of the  Hardy-Littlewood
maximal function was initiated by Juha Kinnunen in \cite{Ki}, where
it was shown that the centered maximal operator is bounded on the
Sobolev spaces $W^{1,p}(\mathbb{R}^{d})$ for $1<p\le\infty$. This
result was used to give a new proof of a weak-type inequality for
the Sobolev capacity, and to obtain the $p$-quasicontinuity of the
maximal function of an element of $W\sp{1,p}(\mathbb{R}\sp{d})$,
$1<p<\infty$. Since then, a good deal of work has been done within
this line of research. In  \cite{KiLi} a local version of the
original boundedness result, valid on $W^{1,p}(\Omega )$,
$\Omega\subset\Bbb R^{d}$ open, was given. Generalizations  were
presented in \cite{HaOn}, extending both the global and local
theorems
 to the spherical maximal function for the range $d>1, d/(d-1)< p$.
 The regularity of
the fractional
maximal operator was studied in \cite{KiSa}.
Hannes Luiro proved in \cite{Lu} the continuity of the Hardy-Littlewood
maximal operator on $W^{1,p}(\Bbb R^{d})$ (continuity is not immediate from
boundedness because of the lack of linearity), finding an explicit representation
for the derivative of the maximal function. Among
other articles  dealing with
related topics we mention \cite{Bu}, \cite{Ko1}, and
\cite{Ko2}.

\

As usual, the case $p=1$ is significantly different from
the case $p>1$, not only because $L^1(\Bbb R^d)$ is not reflexive (so
weak compactness arguments used when $1<p<\infty$ are not available for $p=1$),
but more specifically to this problem, because $Mf\notin L^1(\Bbb R^d)$
whenever $f$ is nontrivial, while the maximal operator acts
boundedly on $L^p$ for $ p>1$.
Nevertheless, in dimension $d=1$, Hitoshi
Tanaka proved (cf. \cite{Ta}) that if
$f\in W^{1,1}(\Bbb R)$, then the noncentered
maximal function  $Mf$ is differentiable a.e. and
$\| DMf\|_1 \le 2\| Df\|_1$. We shall be concerned (mostly) with the case $d= 1$ and $p=1$. What had not previously been  noticed, and we show here,
is that the maximal operator can actually improve the regularity of a function $f$,
rather than simply preserving it, and without increasing the variation. This leads to the possibility of obtaining,
 under less smoothness, versions of classical inequalities involving a function and its derivatives.

\

Recall that if
$f\in W^{1,1}(\Bbb R)$, then $f$ is absolutely continuous and of bounded variation.
We refine Tanaka's
arguments, obtaining the best possible bound
 and generalizing it to the class of functions of bounded variation.
Let $I$ be an interval, let $f: I\to\Bbb R$ be of bounded variation,
and let $Df$ be its distributional derivative.  Denoting by $Mf$ the
noncentered maximal function of $f$, we  prove (cf. Theorem
\ref{main}) that $Mf$ is absolutely continuous. Hence, $Mf$ is
differentiable a.e. and its pointwise derivative coincides with its
distributional derivative $DMf$; thus,  the latter is a function and
not just a Radon measure. Furthermore, the variation of $Mf$ is no
larger than that of $f$, in the sense that $\| DMf\|_{L^1(I)} \le
|Df|( I)$. This inequality is easily seen to be sharp.  Also,
without some assumption of bounded variation type the result fails:
There are bounded, compactly supported (hence integrable) functions
such that $Mf$ is not differentiable on a set of positive measure
(see Example \ref{cantor}).
 We mention that for bounded intervals $I$, the fact that
$\| DMf\|_{L^1(I)} \le |Df|(I)$ tells us that $Mf: BV(I) \to
W^{1,1}(I)$ boundedly.

\

Finally, we note  that from the  viewpoint of regularity the
noncentered maximal operator is better behaved than the centered
one: The latter yields a discontinuous function when applied, for
instance, to the characteristic function of $[0,1]$. And the same
can be said about the one directional maximal operators. In higher
dimensions and for $p > 1$, we mention that even though Kinnunen
stated his boundedness result from \cite{Ki}
 only for the centered operator,
it also holds for the uncentered one by a simple modification of his
arguments, as noted by Tanaka in \cite{Ta}. Alternatively,
boundedness of the uncentered operator can be deduced from Theorem 1
of \cite{HaOn}.

\

This paper is organized as follows. Section 2 contains the basic definitions,
the
main result, and  corollaries. Section 3, the lemmas used in the proof.
In Section 4 examples are presented illustrating the basic issues involved and showing that the main theorem is in some
sense optimal. As an application, in Section 5 we give a variant of Landau's  inequality
under less regularity. Save for the issue of best constants, the inequality we present
is stronger than Landau's. Of course, this kind of argument can be applied to other
inequalities also. As a second, simple instance, we present a trivial variant of the
Poincar\'e-Wirtinger inequality.

\section{Definitions, main theorem and corollaries.}

By an interval $I$ we mean a nondegenerate interval, so examples such as
$[a,a]$ or $(a,b)$ with $b\le a$ are excluded. But we do include the
cases $a=-\infty$ and $b=\infty$, so
$I$ may have infinite length or even be the whole
real line. Let $\lambda$ denote Lebesgue measure and  $\lambda^*$ Lebesgue outer measure.
When $a$ and $b$ are distinct real numbers, not necessarily in increasing order,
we let $I(a,b)$ stand for a (nonempty) interval whose extremes are $a$ and  $b$, while
if the interval $I$ is given, we use $\ell (I)$ and $r(I)$ to denote its left and right
endpoints. Since functions of bounded variation always have lateral limits,
we can go from $(a,b)$ to $[a,b]$ by extension, and viceversa by restriction.
Thus, in what follows it does not really matter whether $I$
is open, closed or neither. Nevertheless, if at some stage of an argument it is useful
to assume that $I(a,b)$ is of a certain type, we shall explicitly say so.

\begin{definition} Given
$P=\{x_1,\dots ,x_L\}\subset I$ with
$x_1 <\dots <x_L$, the variation of   $f:I\to \Bbb R$ associated to the partition $P$ is defined as
$$
V(f, P):= \sum_{j=2}^{L} |f(x_j) - f(x_{j-1})|,
$$
and the variation of $f$ on $I$, as
$$
V(f):=\sup_P V(f,P),
$$
where the supremum is taken over all partitions $P$ of $I$.
We say that $f$ is of {\em bounded variation} if $V(f) <\infty$.
\end{definition}

We use $Df$ to denote the distributional derivative of $f$, and
$I, J$ to denote intervals. Of course, if $f:I\to \Bbb R$ is
absolutely continuous then $Df$ is a function, which coincides with
the pointwise derivative $f^\prime$ of $f$. In this case we also
denote the latter by
$Df$.

\begin{definition} The canonical representative of $f$ is the
function
$$
\overline{f}(x) := \limsup_{\lambda (I)\to 0, x\in I}\frac{1}{\lambda (I)}\int_I f(y)dy.
$$
\end{definition}

The Lebesgue differentiation
theorem tells us that  $\overline{f}=f$ a.e., so $\overline{f}$
does represent the equivalence class of $f$. Of course, taking the
$\liminf$ would yield a representative as ``canonical" as the one
above; we just selected the one best suited to our purposes.

It is well known (cf. Lemma \ref{equivalence} in the next section)
that if $f$ is of bounded variation, then $Df$ is a Radon measure
with $|Df|(I) = V(\overline{f}) <\infty$, where $|Df|$ denotes the
total variation of $Df$.

\begin{definition} Given  a locally integrable function $f:I\to \Bbb R$, the
noncentered Hardy-Littlewood
maximal function $Mf$  is defined by
$$
Mf(x) := \sup_{ x\in J\subset I}\frac{1}{\lambda (J)}\int_J
|f(y)|dy.
$$
If $R > 0$, the definition of the {\em local} maximal function $M_Rf$ is the same as above, save
that the intervals $J$ are also required to satisfy $\lambda J \le R$.
\end{definition}

\begin{remark}  The terms {\em restricted} and {\em truncated}
have been used in the literature to designate $M_Rf$. However, in both cases the meaning differs
from the usual notions of restriction and truncation of a function, so we
prefer {\em local}. $M_Rf$ is genuinely local in that its value at $x$ depends
only on how $f$ behaves in an $R$-neighborhood of $x$.

For some purposes the relevant maximal operator is $M_R$ rather than
$M$ (cf., for instance, \cite{Ha}). Thus, it seems worthwhile to point out
that the results we prove on regularity and size of the derivative
 hold both for $M$ and $M_R$, essentially for the same reasons.
 So in the
statements of the theorems $M$ and $M_R$ will appear, but the proofs will
only mention $M$, unless some modification is needed to cover the
case of $M_R$. With respect to the possibility of deriving the results for $M_R$ from
those of $M$, or viceversa,
it is not clear to us whether this can be done, in view
of the fact that
 neither $DMf$ nor $ DM_Rf$ pointwise dominates  the other.

The local maximal
function will be used at the end of
this paper to prove an inequality of Poincar\'e-Wirtinger type. There $M_R$ cannot
be replaced by $M$.
\end{remark}

\begin{theorem}\label{main}  If $f:I\to \Bbb R$ is of
bounded variation
then
 $Mf$ is absolutely continuous. Furthermore,
$ V(Mf) \le V(\overline{f}), $ or equivalently, $ \|DMf\|_{L^1(I)}
\le |Df|(I).$ The same holds for $M_Rf$.
\end{theorem}
\begin{proof} We  assume that
$0\le f =\overline{f}$. Since $f$ is upper semicontinuous (Lemma
\ref{semicont}), and $f\le Mf$, the maximal function
 is continuous (Lemma \ref{Mcont} and Remark \ref{local}), and of
bounded variation with $V(Mf)\le V(f)$ (Lemma \ref{vardec}). Also,
the image under $Mf$ of a measure zero set has measure zero (Lemma
\ref{nullimage}), so by the Banach Zarecki Theorem (Lemma
\ref{BVAC}) $Mf$ is absolutely continuous, whence $|DMf|(I)=
\|DMf\|_{L^1(I)}$. Finally, by Lemma
 \ref{equivalence}, $|DMf|(I) = V(Mf)$ and $|Df|(I) =  V(f)$, so
 $\| DMf\|_{L^1(I)} \le |Df|(I)$.
 \end{proof}

\begin{remark} Actually, the proof of Lemma \ref{vardec} yields a slightly stronger
result: Given $f$ and $Mf$ on $I$, and any subinterval $J\subset I$
such that the endpoints of $J$ belong to $\{Mf=f\}$, we have
$V(Mf|_J)\le V(f|_J)$. That is, the variation fails to grow not only
when considered over the whole interval, but also over a wide class
of subintervals. The key to this reduction of the variation is  the
rather simple behaviour of $Mf$ on the components of the open set
$\{Mf
> f\}$: If $(\alpha, \beta)$ is any such component, then either $Mf$ is
monotone there, or there exists a $c$ in $(\alpha, \beta)$ such that
$Mf$ is decreasing on $(\alpha, c)$ and increasing on $(c, \beta)$.
\end{remark}

We recall  the definitions of the space $BV(I)$ of
integrable functions of bounded variation and of the
Sobolev space $W^{1,1}(I)$.

\begin{definition} Given the interval $I$,
$$
BV(I) := \{f:I\to \Bbb R| f\in L^1(I) \mbox{ and } |Df|(I) <\infty\},
$$
and
$$
W^{1,1}(I) := \{f:I\to \Bbb R| f\in L^1(I), Df \mbox{ is a function, and } Df\in L^1(I)\}.
$$
\end{definition}

It is obvious that $W^{1,1}(I) \subset BV(I)$ and that the inclusion
is proper. The Banach space $BV(I)$ is endowed with the norm
$\|f\|_{BV(I)}:= \|f\|_{L^1(I)} + |Df|(I)$, and $ W^{1,1}(I),$ with
the restriction of the $BV$ norm, i.e.,  $\|f\|_{W^{1,1}(I)}:=
\|f\|_{L^1(I)} + \|Df\|_{L^1(I)}$.

\begin{remark} It is well known that for every $f\in L^1(I)$ and every $c > 0$,
$c \lambda (\{Mf>c\}) \le 2 \|f\|_{L^1(I)}$, so by Theorem
\ref{main}, the maximal operator satisfies a ``mixed type"
inequality on $BV(I)$, weak type for functions and strong type for
their derivatives. Hence $M$ maps $BV(I)$ into the subspace of
$L^{1,\infty}$ consisting of the functions whose distributional
derivative is an integrable function. But if $I$ is bounded,
something stronger can be said: $Mf$ maps boundedly $BV(I)$ into
$W^{1,1}(I)$.
\end{remark}

\begin{corollary}\label{BVtoSob2} Let $I$ be bounded. Then there exists a constant $c= c(I)$ such that
for every $f\in BV(I)$,  $Mf\in W^{1,1}(I)$ and $\|Mf\|_{W^{1,1}(I)}\le  c \|f\|_{BV(I)}.$
\end{corollary}
\begin{proof} By Sobolev embedding
for $BV$ functions (cf. Corollary 3.49 p. 152 of \cite{AFP}),
$\|f\|_\infty\le  c(I) \|f\|_{BV(I)}$, so
$$
\|Mf\|_{L^1(I)}\le \lambda(I) \|f\|_\infty\le \lambda(I) c(I)
\|f\|_{BV(I)}.
$$\end{proof}

After this article was completed, the authors were able to show (cf.
\cite{AlPe}, Theorem 2.7) that the local maximal operator $M_R$ is
bounded from $BV(I)$ into $W^{1,1}(I)$ {\it even} if $I$ has
infinite length; in addition, the bounds  grow with $R$ as $O(\log
R)$.

\section{Lemmas.}

 Given
$f:I\to \Bbb R$, define the
upper derivative of $f$ as
$$
\overline{D}f(x)= \limsup_{h\to 0}\frac{f(x+h) - f(x)}{h}.
$$

The next lemma is well known and we include it here for the reader's convenience.

\begin{lemma}\label{image}
 Let $f:I\to \Bbb R$ be a continuous function, and let
$E\subset \left\{x\in I: \left|\overline{D}f(x)\right| \le k\right\}$.
Then
 $\lambda^*\left( f\left(E\right)\right) \le k\lambda^* E$.
\end{lemma}
\begin{proof}
We may assume that I is open. Fix $\varepsilon > 0$, and let
$O\subset I$ be an open set with $E\subset O$ and $\lambda O \le
\lambda^* E +\varepsilon$. For each $x\in E$, pick $y$ so that the
closed interval $I(x,y)$ is contained in $O$ and for
 every $z\in I(x,y)$,
 $\left|\frac{f(z) - f(x)}{z-x}\right| \le k + \varepsilon $.
Then ${\mathcal V} :=\{I(x, z): x\in E, z\in I(x,y) \subset O\}$
is a covering of $E$ with $\lambda ( \cup{\mathcal V})\le \lambda O
\le \lambda^* E +\varepsilon$.
By continuity of $f$, for every pair $\{x,z\}$ the set $f(I(x,z))$  is connected,
hence an interval or a point. Let $C :=\{c\in {\Bbb R}: \mbox{ for some } I(x,z)\in {\mathcal V}, f(I(x,z))= c\}$.
Suppose  $c_1, c_2\in C$ and $c_1\ne c_2$. Then $f^{-1} (\{c_1\})$ and $f^{-1} (\{c_2\})$
are disjoint sets and each contains an interval, from which it follows that $C$ is at most countable.
Now if $z\to x$ then
$f(z)\to f(x)$, so the collection  $f({\mathcal V^\prime}) :=\{f(I(x, z)):
 I(x,z) \in{\mathcal V} \mbox{ and } f \mbox{ is not constant on } I(x,z)\}$ is a Vitali covering of
$f(E) \setminus C$. Hence there is a
 disjoint subcollection   $\{I_n\}$ of $f({\mathcal V^\prime})$ such that
$\lambda^*\left(f( E)\setminus \cup_n I_{n})\right) = 0$.
For each $n$, select $I(x_n, z_n)$ such that $f(I(x_n, y_n)) = I_n$.
Then
$$\lambda^*\left(f( E)\right) \le \sum_n \lambda (I_{n}) \le \sum_n (k+\varepsilon ) \lambda (I(x_{n}, y_n))\le   (k +\varepsilon) \lambda O
\le ( k +\varepsilon) (\lambda^*E +\varepsilon) .$$
\end{proof}

The following lemma  is the direction we need of the Banach Zarecki Theorem (an
``if and only if"  result). It is an immediate consequence of the Fundamental Theorem
of Calculus for the Lebesgue integral, the a.e. differentiability of functions
of bounded variation, and the preceding lemma.

\begin{lemma}\label{BVAC}
 Let $f:I\to \Bbb R$ be a continuous function of bounded
variation. If $f$  maps measure zero sets
into measure zero sets, then it is absolutely continuous.
\end{lemma}

Let $f:I\to \Bbb R$ be a
function of bounded variation. Then
$|f|$ is also of bounded variation, and
$\left|D|f|\right|( I ) \le |Df|( I )$.
Additionally,
when studying the boundedness properties of the
maximal function it makes no difference whether we consider
$f$ or $|f|$. For notational simplicity we will often
assume that
$f\ge 0$.

The next two lemmas are likely to be well  known, and small variants
certainly are. But since we
are not aware of any explicit written reference, we include
them for completeness.

\begin{lemma}\label{semicont}
The canonical representative $\overline{f}$ of a function of bounded variation $f:I\to \Bbb R$ is upper semicontinuous and of bounded variation. Furthermore,
$\overline{f}$ minimizes the variation within the equivalence class of $f$.

\end{lemma}

\begin{proof}
Recall that a function $f$ of bounded variation (being the difference of two monotone functions)
has left and right limits at every point.
To see that $\overline{f}$
is upper semicontinuous,
simply note that
$\overline{f}(x)=\max \{\lim_{y\uparrow x}\overline{f}(y), \lim_{y\downarrow x}\overline{f}(y)\}$, so for every $x\in I$ and
every sequence $\{x_n\}$ in $I$ converging to $x$, $\limsup_n\overline{f}(x_n)\le \overline{f} (x)$.
And $V(\overline{f})\le V(f)$ follows immediately
from the fact that $\overline{f} (x)$ satisfies the following two conditions:
i) If $x$ is a
point of continuity of $f$, then $\overline{f} (x) = f(x)$; ii)  if $x$ is a point of discontinuity of $f$, then
$\overline{f}(x)$ belongs to the closed interval determined by the extremes $\lim_{y\uparrow x}f(y)$ and  $\lim_{y\downarrow x}f(y)$. Finally, it is quite obvious (or else, cf. Theorem 3.28,  page 136 of \cite{AFP})
that on the equivalence class of $f$,
$V$ achieves its minimum value at $g$ iff  $g$ satisfies both conditions i) and ii) above.
\end{proof}

Next, we consider balls associated to
some norm in $\Bbb R^d$ and the corresponding maximal operator.
In the local case $R$ will denote the diameter, in accordance with our one
dimensional convention. While in this paper we only need the case $d=1$,  we state the next result for arbitrary $d$, as it
 is likely to be useful in future work.

\begin{lemma}\label{Mcont}
 Let $f:\mathbb{R}^d\to [0,\infty]$ be locally integrable.
 If $f$ is upper semicontinuous at
 $w\in \mathbb{R}^d$ and $f(w)\le Mf(w)$, then $Mf$ is continuous at $w$.
The same holds for $M_Rf$.
\end{lemma}
\begin{proof}
Since $Mf$ is lower semicontinuous, it suffices to prove the upper semicontinuity
of $Mf$ at $w$. The idea is simply to note that if  ``large" balls are considered near
$w$, increasing their radii a little so as to include $w$ cannot decrease
the average by much, while if one is forced to consider arbitrarily
small balls,  then the fact that $\limsup_n f(x_n) \le f(w)$
whenever $x_n\to w$   leads to the same
conclusion for $Mf$. More precisely, we show that
given $\varepsilon >0$, there exists a $\delta >0$ such
that for all $x \in B(w,\delta)$, $Mf(x) \le Mf(w)+\varepsilon$.
Fix $k>>1$ such that $(1+1/k)^d Mf(w) \le Mf(w) +\varepsilon$, and choose $\tau >0$ with
$f(y)\le f(w)+\varepsilon$ whenever $|y -w| <3\tau$. Set $\delta = \tau /k$,
let $|x - w|<\delta$, and let $B(u,r)$ be a ball containing $x$.
If $B(u,r)\subset \{f\le f(w) + \varepsilon\}$, then
\begin{equation*}
\frac{1}{\lambda (B(u,r))}\int_{B(u,r)} f(y)dy \le f(w)+\varepsilon \le Mf(w)+\varepsilon,
\end{equation*}
while if $B(u,r)\cap \{f\le f(w) + \varepsilon\}^c \ne \emptyset$,
then $r>\tau$, and since
$w\in B(u,r+\delta )$, we have
$$
\frac{1}{\lambda (B(u,r))}\int_{B(u,r)} f(y)dy
\le
\frac{\lambda (B(u, r+\delta ))}{\lambda (B(u,r))}
\frac{1}{\lambda (B(u,r+\delta ))}\int_{B(u,r+\delta )} f(y)dy
$$
$$
\le
\left(1 + \frac{\delta}{r}\right)^d  Mf(w)
\le Mf(w) +\varepsilon.
$$
In the case of $M_Rf$, if the large balls already have diameter $R$, instead
of increasing the radii just translate the balls slightly. The easy details
are omitted.
\end{proof}

\begin{remark}\label{local} Actually, we shall apply the preceding
lemma to arbitrary intervals $I\subset \Bbb R$ and not just to $\Bbb R$.
Nevertheless, the result is stated for $\Bbb R^d$ since when $d > 1$, there are
connected and simply connected open sets $O\subset\Bbb R^d$ for which it fails, as we
shall see. The difference stems from the fact that when working on $O$, the maximal
function is defined by taking the supremum over balls contained in $O$.  If $d = 1$,
the proof of the lemma can be easily adapted to intervals. Alternatively, the result
for $\Bbb R$ implies the general case as follows: Given  $I\subset \Bbb R$, extend
$f\ge 0$ to $\Bbb R$ by setting it equal to zero off $I$. This changes neither the
upper semicontinuity of $f$ at points in $I$, nor the values of $Mf$ on $I$. Regarding
 the case $d > 1$, counterexamples  already exist when
$d = 2$. For convenience we use the $\ell_\infty$ norm $\|(x,y)\|_\infty =
\max\{|x|,|y|\}$ on
$\Bbb R^2$, so balls will refer to the metric defined by $\|\cdot\|_\infty$. Let
$O\subset
\Bbb R^2$ be the open set $(0,2)^2 \cup (0,3)\times (0,1)$. Define $f:=
\chi_{(0,1]^2}$ on $O$. Then $f$ is upper
semicontinuous, $Mf\ge 1/4$ on
$(0,2)^2$, $Mf = 0$ on $(2,3)\times (0,1)$ (so $Mf$ is not continuous), and $Mf\ge
f$ everywhere.
\end{remark}

\begin{lemma}\label{locmax}
Let  $f:I\to [0,\infty)$ be  an upper semicontinuous function such
that for every $x\in I$, $f(x)\le Mf(x)$. Suppose there exists an
interval $[a,b]\subset I$ with $\max \{Mf(x): x\in [a,b]\}  >
\max\{Mf(a), Mf(b)\}$. If $c\in [a,b]$ is a point where $Mf$
achieves its maximum value on $[a,b]$, then $Mf(c) = f(c) =
\max_{x\in [a,b]} f(x)$. If either $f\in L^1(I)$ or we consider
$M_Rf$ instead of $Mf$, then the same result holds under the
following weaker assumption: $Mf$ (respectively  $M_Rf$) achieves
its  maximum value on $[a,b]$ at some interior point, so  $\max
\{Mf(x): x\in (a,b)\}  \ge \max\{Mf(a), Mf(b)\}$ (respectively $\max
\{M_Rf(x): x\in (a,b)\}  \ge \max\{M_Rf(a), M_Rf(b)\}$).
\end{lemma}

\begin{proof} By Lemma \ref{Mcont}, $Mf$ is continuous.
Suppose that $Mf(c) =  \max \{Mf(x): x\in [a,b]\} > \max\{Mf(a), Mf(b)\}$. If
$f(c) < Mf(c)$, by upper semicontinuity of $f$ there exists an open interval
$J:=(c-\delta, c+\delta)$ such that if $x\in J$, then $f(x)<2^{-1}(f(c) +
Mf(c))$.  Define $L(t) = \frac{1}{t}\int_{c-t}^c f$ on $[\delta ,c-a]$ and
$R(x) = \frac{1}{x}\int_{c}^{c+x} f$ on $[\delta, b-c ]$. By continuity, there
exist
$t_0$ and $x_0$ maximizing $L$ and $R$ respectively.
Since $Mf(c) > \max\{Mf(a), Mf(b)\}$,
in order to evaluate $Mf(c)$ we only need to
consider intervals properly contained in $[a,b]$, so
$Mf(c) = \max\{L(t_0), R(x_0)\}$. Suppose without
loss of generality that $L(t_0)\le R(x_0)$.
Then
$$
Mf(c) = \frac{1}{x_0}\int_{c}^{c+x_0} f \le
\frac{1}{x_0}\left(\int_{c}^{c+\delta}2^{-1}(f(c) + Mf(c)) +\int_{c+ \delta}^{c+x_0} Mf(c)\right) < Mf(c),
$$
 a contradiction.

Suppose next that  $M_Tf$ achieves its maximum value on $[a,b]$
at some interior point, so $\max \{M_Tf(x): x\in (a,b)\}  \ge \max\{M_Tf(a),
M_Tf(b)\}$.
Define $L$ on $[0,\min\{T, c-\ell (I)\}]$ and $R$ on
$[0, \min\{T, c+r (I)\}]$ as above, and conclude that there exist $t_0$ and
$x_0$ maximizing $L$ and
$R$ respectively. Then argue as before. If  $f\in L^1(I)$ we reason in the same
way. Suppose for instance  that $I=\Bbb R$. Then
$\lim_{t\to \infty}L(t) =\lim_{x\to \infty}R(x) = 0$, so again there are points $t_0$ and $x_0$ maximizing $L$
and $R$.  The case of bounded or semi-infinite intervals is easily
handled, assuming, for instance, that finite extremes belong to $I$, and concluding
as before that $t_0$ and $x_0$ exist.
 \end{proof}

 \begin{lemma}\label{equivalence}  Let $f:I\to [0,\infty)$ be  of bounded variation.
 Then $V(\overline{f}) = |Df|(I)$.

\end{lemma}
\begin{proof} This is well known, and it follows from \cite{AFP}, Proposition 3.6 p. 120 together
with Theorem
3.27 p.135 (making the obvious adjustments in the definitions, to take into
account that we do not
assume  $f\in L^1(I)$).
\end{proof}

In what follows we will distinguish between the points where one needs to consider arbitrarily
small intervals to obtain the value of the maximal function, and those where the supremum is achieved
by looking at intervals of length bounded below. Since at this stage we are not
assuming that $f$ is integrable, it might happen that for some increasing sequence of
intervals $I_n$ containing $x$, and for instance, with
$I_n\uparrow I = [a,\infty)$,
$Mf(x) >  \frac{1}{\lambda (I_n)}\int_{I_n} f$ but
$Mf(x) = \sup_n  \frac{1}{\lambda (I_n)}\int_{I_n} f$.
As a shorthand to describe this situation, we write
$Mf(x) = \frac{1}{\lambda (I)}\int_{I} f$, and similarly for
other intervals of infinite length. Now
set
$$E:= \left\{x\in \Bbb R:\mbox{ there exists an interval } I \mbox{ containing } x
\mbox{ such that } Mf (x)= \frac{1}{\lambda (I)}\int_I f\right\}.
$$
In the remaining of this section
we assume once and for all that $f = \overline{f}$ whenever $f$ is of bounded variation.
Then  the Lebesgue theorem on differentiation of integrals together with
 $0\le f = \overline{f}$,
entail that $E^c\subset \{Mf = f\}$. Next, write
$$E_n:= \left\{x\in \Bbb R:  \mbox{ there exists an } I
\mbox{ such that } x\in I,  Mf (x)= \frac{1}{\lambda (I)}\int_I f, \mbox{ and }
\lambda (I)\ge\frac{1}{n}\right\}.$$

If we are dealing with  $M_Rf$, the sets $E$ and $E_n$ are defined as before save for the fact that we add the extra condition $\lambda (I) \le R$ (so for instance, if $R=1/2$, $E_1=\emptyset$).  As usual, Lip$(g)$ denotes the Lipschitz constant of a Lipschitz function $g$.

\begin{lemma}\label{Lip}
Let $f\in L^\infty (I)$. Then the restriction of $Mf$  to $E_n$ is  Lipschitz, with
{\rm Lip}$(Mf) \le n \|f\|_\infty$. The same holds for $M_Rf$.
\end{lemma}
\begin{proof} Fix $x, y\in E_n$. By symmetry, we may assume that $Mf(x) \ge Mf(y)$.
 By hypothesis, there exists a $J$   containing $x$
 such that $\lambda (J)\ge\frac{1}{n}$ and $Mf (x)= \frac{1}{\lambda (J)}\int_J f$.
Our notational convention allows $J$ to have infinite length, so we suppose first
that $\lambda (J)<\infty$. Then
$$
\frac{Mf(x) - Mf(y)}{|x-y|}\le \frac{\frac{1}{\lambda (J)}\int_J f -
\frac{1}{\lambda (J) + |x-y|}\int_J f}{|x-y|}
= \frac{\frac{1}{\lambda (J)}\int_J f}{\lambda (J) + |x-y|}
\le \frac{\|f\|_\infty}{\lambda (J) +|x-y|}< n \|f\|_\infty.
$$
The case $\lambda (J)=\infty$ is obtained by an easy approximation argument.

With respect to $M_R$, the reasoning is similar. Suppose
that $x, y\in E_n$, with $x<y$.
If $y -x \ge 1/n$, then
$$
\frac{|M_Rf(x) - M_Rf(y)|}{|x-y|}
\le \frac{\|f\|_\infty-0}{y-x}\le  n \|f\|_\infty.
$$
So assume that $y -x < 1/n$ and $M_Rf(x) > M_Rf(y)$ (the case $M_Rf(x) < M_Rf(y)$
is
handled in the same way). By hypothesis there exists an interval $[a,b]$ containing
$x$ such that
$1/n \le b-a \le R$ and
 $$
 M_Rf (x)= \frac{1}{b-a}\int_a^b f.
 $$
Now if $y - a \le R$ we can repeat the argument given for $Mf$, so suppose $y - a >R$.
Writing $c := R - b +a$, we have
$$
\frac{M_Rf(x) - M_Rf(y)}{y-x}\le \frac{\frac{1}{b-a}\int_a^b f -
\frac{1}{R}\int_{y-R}^y f}{y-x}
\le \frac{(b-a)\int_a^b f - (b-a)\int_{y-R}^b f + c \int_a^b f}{(y-x)(b-a) R}
$$
$$
\le \frac{(y-R-a)\|f\|_\infty   + c \|f\|_\infty}{(y-x) R}
= \frac{y-b}{y-x}\frac{\|f\|_\infty}{R}\le n \|f\|_\infty.
$$

\end{proof}

\begin{lemma}\label{vardec}
Let  $f:I\to [0,\infty)$ be  an upper semicontinuous function such that for every $x\in I$, $f(x)\le Mf(x)$. Then
$V(Mf)\le V(f)$. In particular,  $Mf$ is of bounded variation whenever $f$ is.
The same results hold for $M_Rf$.
\end{lemma}
\begin{proof}  We show that $Mf$ varies no more  than
$f$ on $\{Mf>f\}$, and of course the same happens on $\{Mf>f\}^c = \{Mf=f\}$.
Note that $\{Mf>f\}$ is open in $I$, since $Mf - f$ is lower semicontinuous, being the
difference between a continuous and an upper semicontinuous function.
Let $(\alpha,\beta )$ be any component of $\{Mf>f\}$
(here we can either assume directly that $I$ is open, or else consider
 as possible components  intervals of the form $[\alpha,\beta )$ and $(\alpha,\beta ]$, which can be handled in the same way as we do below). We will see next
that on
$(\alpha,\beta)$, $Mf$ can
only behave  in one of two ways:
Either $Mf$ is monotone, or
there exists a $c \in (\alpha,\beta)$ such that $Mf$
decreases on $(\alpha,c)$ and increases on $(c,\beta)$.
Suppose towards a contradiction that for some points $c_1, c_2, c_3$
with $\alpha < c_1< c_2< c_3 < \beta$ we have
$Mf(c_1)< Mf(c_2)$ and $Mf(c_{3})< Mf(c_2)$.
By changing $c_2$ if needed, and relabeling,
we may assume that $Mf(c_2) =\max\{Mf(x): x\in [c_{1},c_3]\}.$ Then $c_2\in \{Mf>f\}$,
and
$
Mf(c_2) = f(c_2)
$ by Lemma \ref{locmax}.
The result for $M_Rf$  also follows from  Lemma
\ref{locmax}.

To show that $V(Mf)\le V(f)$, given  an arbitrary partition
$\{x_1,\dots ,x_L\}$ of  $I$ we produce a refinement
$\{y_1,\dots ,y_K\}$ such that
$\sum_2^K |Mf(y_i) - Mf(y_{i-1})|
\le \sum_2^K |f(y_i) - f(y_{i-1})|$.
We always assume that partitions are labeled in increasing order.
Of course, if $\{x_1,\dots ,x_L\}\subset \{Mf= f\}$ there is nothing to do.
Otherwise, for each $x_i\in \{Mf > f\}$ there is a unique component
$(\alpha_i , \beta_i )$ which contains it. Add both endpoints
$\alpha_i$ and $\beta_i$ to the partition.
If $Mf$ is monotone on
$(\alpha_i , \beta_i )$, we do not add any new points inside this interval.
If $Mf$ is {\em not} monotone on $(\alpha_i , \beta_i )$,
then there exists a $c_i \in (\alpha_i,\beta_i)$ such that $Mf$
decreases on $(\alpha_i,c_i)$, increases on $(c_i,\beta_i)$,
and $Mf(c_i) < \min \{Mf(\alpha_i) ,Mf(\beta_i)\}$. If $c_i$ does
not already belong to the original partition, include it in the refinement. The resulting finite collection
$P:=\{y_1,\dots ,y_K\}$  has the desired
property, as we shall see. First we study what
happens in each $[\alpha_i , \beta_i ]$.
For every pair $y_{i} < y_{i+k}$
of points of $P$ that are endpoints
of some component $(\alpha_{j_i} , \beta_{j_i} )$ and contain
one or more elements of $P$ between them, say
$y_{i} < y_{i+1}<\dots  < y_{i+k}$, either $Mf$ is monotone
on that component, and then
\begin{equation*}
\sum_{j=i+1}^{i+k} |Mf(y_j) - Mf(y_{j-1})|
= |Mf(y_{i+k}) - Mf(y_{i})|
\end{equation*}
\begin{equation*}
 = |f(y_{i+k}) - f(y_{i})|
\le \sum_{j=i+1}^{i+k} |f(y_j) - f(y_{j-1})|,
\end{equation*}
or $Mf$ achieves its minimum value on $[y_{i}, y_{i+k}]$ at some
intermediate $y_{i_m}$ and $Mf(y_{i_m}) < \min \{Mf(y_{i})
,Mf(y_{i+k})\}$. In this case
\begin{equation*}
\sum_{j=i+1}^{i+k} |Mf(y_j) - Mf(y_{j-1})|
= |Mf(y_{i+k}) - Mf(y_{i_m})| + |Mf(y_{i_m}) - Mf(y_{i})|
\end{equation*}
\begin{equation*}
< |f(y_{i+k}) - f(y_{i_m})| + |f(y_{i_m}) - f(y_{i})|
\le \sum_{j=i+1}^{i+k} |f(y_j) - f(y_{j-1})|.
\end{equation*}
Finally, for each pair $\{y_i, y_{i+1}\}$
not already taken into account, we have
$|Mf(y_{i+1}) -Mf(y_i)| = |f(y_{i+1}) -f(y_i)|$,  so
the conclusion follows.
\end{proof}

Before proving the next lemma, we mention that on large
parts of its domain $Mf$ is  locally Lipschitz.
Of course, matters would be considerably simpler if
$Mf$ were locally Lipschitz at every
point, but
unfortunately this need not be the case, as the following
example shows: Take $f(x) = (1-\sqrt x )\chi_{[0,1]}
(x)$ and note that
$Mf$ fails to be locally Lipschitz at $0$.

\begin{lemma}\label{nullimage}
 Let $f:I\to \Bbb R$ be a  function of bounded variation, and let
$N$ be a set of measure zero.
Then
 $\lambda\left(M f\left(N\right)\right) =0$. The same results hold for $M_Rf$.
\end{lemma}
\begin{proof}  Suppose
$N$ has measure zero, and let
 $E_n$ and $E$ be the sets whose
definition appears just before Lemma \ref{Lip}.
Since $E_n\uparrow E$, by Lemmas \ref{image} and
\ref{Lip},
for each $n = 1,2,\dots$,
$\lambda (Mf(N\cap E_n)) = 0$, so
$\lambda (Mf(N\cap E)) = 0$. Thus,
 we  may assume that $N\subset E^c$, whence
$N\subset \{Mf = f\}$. We are going to make further reductions
on $N$. First, we remove from it all the intervals
$I_\alpha$ where $Mf$ is constant; we
can do that since $\lambda Mf(\cup_\alpha I_\alpha) =0$
by Lemma \ref{image} (note that the
difference between considering closed or open intervals is
at most a countable set of endpoints, so the exact nature of the
$I_\alpha$'s is of no consequence here). And second, we eliminate from
$N\setminus  \cup_\alpha I_\alpha$ a countable
set in such a way that every remaining point is a
point of accumulation. This can be done
by a well known argument: Pick a countable  base ${\mathcal B}$ of intervals, and
let $\{I_j\}$ be the collection of all intervals in
${\mathcal B}$ for which $(N\setminus \cup_\alpha I_\alpha)\cap I_j$
is countable. Then $(N\setminus  \cup_\alpha I_\alpha)\setminus \cup_jI_j$
has the desired property. For the usual reason of
notational simplicity, we use $N$ again to denote this
thinner set.

Now we are ready to suppose, towards a contradiction, that
 $\lambda^* Mf (N)  > 0$.
Write $8c:= \lambda^* Mf (N)$.
We show that for every finite sequence of distinct real numbers
$\{x_1,\dots ,x_L\}$ in $I$, labeled in (strictly) increasing order,
there is a refinement
$\{y_1,\dots ,y_K\}$ in $I$ with $y_1 < y_2 <\dots <\ y_K$ and
$\sum_2^L |f(x_i) - f(x_{i-1})| + c < \sum_2^K |f(y_i) - f(y_{i-1})|$.
This contradicts the fact that  $f$
 is of bounded variation.

 The final partition will be produced in several
stages, so at any step in the argument, we use
$P$ to denote the partition already at hand, and
$P^\prime$ the immediate refinement obtained in that
step. So $P$ will denote different partitions
at different stages, and the same happens with
$P^\prime$.

By adding more points if needed, and relabeling in increasing order,
 we may assume that
``a large part" of $N$ is contained between the first and
the last points of the partition (call them
$A$ and $B$ respectively). By this we mean that
 $7 c< \lambda^* Mf (N\cap [A,B])$. Again we use $N$ to denote the
 null set $N\cap (A,B)$, and  $\{x_1,\dots ,x_L\}$
 to denote the points of the new partition.

When we say that an  interval
$J$ is determined by $P$ we mean that its extremes are consecutive
points in $P$.
Let  $I_1$ be the first (according to the real ordering)
 of the open intervals determined by
$P:=\{x_1,\dots ,x_L\}$ such that
$\lambda^* Mf (N\cap I_1) > 0$, and let $x_{i_1} < x_{i_1 + 1}$
be the  endpoints of this interval.
For each $x\in N\cap I_1$ and $n\in \Bbb N$ pick $y_n\in N\cap I_1$
such that $\lim_n y_n =x$ and
$|x-y_n| < 2^{-1}  d(x_{i_1+1},I(x, y_n))$
(where $d$ stands for distance, and $I(x, y_n)$ for the {\em compact} interval with extremes $x$ and $y_n$). In the case of $M_Rf$ we additionally require that for every $x$ and every
$y_n$ in its associated sequence $\{y_n\}_{n=0}^\infty$, $2 |x-y_n| < R$.
Since
$Mf$ is continuous and
 nonconstant in all of those
intervals, $Mf(I(x, y_n))$ is  a nondegenerate compact
interval and  $\lim_n \lambda
\left(Mf(I(x, y_n))\right) = 0$.
Thus the collection  ${\mathcal V} :=\{Mf(I(x, y_n)): x\in N, n\in{\mathbb N} \}$
is a Vitali covering of
$Mf(N\cap I_1)$. Furthermore, if $J$ is the closed interval
with the same left endpoint as $I(x, y_n)$ and twice its length,
then $J\subset I_1$. Select
a finite, disjoint subcollection   $\{S_1,\dots ,S_{R}\}$ from
${\mathcal V}$,  such that
\begin{equation}
\label{size}
\lambda^*\left(Mf(N\cap I_1)
\setminus (S_{1}\cup\dots \cup S_{R})\right) < 7^{-1}\lambda^* Mf(N\cap I_1).
\end{equation}
For
each $i= 1, \dots , R$ pick
$J_i\in \{I(x, y_n): x\in N\cap I_1, n\in{\mathbb N} \}$ with $Mf(J_i) = S_i$.

Without loss of generality, suppose that $f(x_{i_1})\le f(x_{i_1+1})$.
Adding a finite number of points between
$x_{i_1}$ and $x_{i_1+1}$ to the original partition $P$
 does not increase the variation  if $f$ is
behaving monotonically there.
Thus, our strategy consists in selecting new points so that
``broken line configurations" are obtained sufficiently often.

We consider two cases. In the first, the increase in the variation is obtained by
adding to the original partition the endpoints of suitable intervals, and either one
or two points inside each such interval. In the second, we add the endpoints of
 intervals not considered in case 1, together with  either one
point {\em outside} each such interval (but close to it), or no additional point.

 Note that $f$ and $Mf$ take the same values
on the endpoints of the intervals $J_i$, since  $\ell (J_i), r (J_i)\in N$
for every $i= 1, \dots , R$.

{\em Case 1.}
Call an interval $J_i$ of type A if $|Mf(r(J_{i})) - Mf(\ell(J_{i}))| >
2^{-1} \lambda S_{i}$ and there exists a
$c_i\in J_{i}$ with $f(c_i) < \min\{f(\ell
(J_{i})), f(r (J_{i}))\} -  2^{-1}|f(r (J_{i})) - f(\ell (J_{i}))|$.
Note that if we have any partition $P$ with consecutive
points $a<b$ and $J_i\subset (a,b)$,
then adding
$\ell (J_{i}), c_i$ and $r (J_{i})$ to $P$ leads to $V(f,P)  +
\frac12 \lambda Mf(J_i)  < V(f, P^\prime)$
(where $P^\prime$ is the refinement so obtained),
regardless of the values of $f(a)$  and $f(b)$.

We say that $J_i$ is of type B if
$|Mf(r(J_{i})) - Mf(\ell(J_{i}))| \le 2^{-1} \lambda S_{i}$.
Suppose in this case that $Mf|_{J_{i}}$ achieves its extreme values at
 $m_1,m_2\in J_{i}$, with $m_1<m_2$. We add the distinct elements in
$\{\ell (J_{i}), m_1, m_2,r(J_{i})\}$ to
$P$,
obtaining $P^\prime$ (it may happen
that either $\ell (J_{i})= m_1$ or $m_2= r(J_{i})$,
but not both, so $P^\prime$ contains either
$3$ or $4$ points more than $P$). Note that
$\lambda S_{i}=|Mf(m_{2}) - Mf(m_{1})| \le
|f(m_{2}) - f(m_{1})|$, since
where the maximum of $Mf$ occurs,
$Mf$ and $f$ take the same value
(by Lemma \ref{locmax} if the corresponding $m_i$ is an interior point,
and by $N\subset \{Mf=f\}$ otherwise) while $f\le Mf$ always.
As before, for some pair of consecutive
points $a<b$ in $P$ we have $J_i\subset (a,b)$.
It is again clear that no matter what the positions of
$f(\ell (J_{i})), f(m_1), f(m_2),$ and $f(r(J_{i}))$ are
relative to $f(a)$ and  $f(b)$,
we always  have
$V(f,P)  +
\frac12 \lambda Mf(J_i)  \le V(f, P^\prime)$.

Suppose now that the collection of intervals $J_{i_j}$ of type either A or B
satisfies $\sum_j \lambda S_{i_j} \ge 3^{-1} \sum_{i=1}^R\lambda S_{i}$.
Adding to the initial partition all their endpoints
and the interior points corresponding
to each case we get
$$
V(f,P) +
\frac{1}{6}\sum^R_1 \lambda S_i \le
V(f,P^\prime).
$$

{\em Case 2.}
If the case previously considered does not
hold, then  the set of intervals $J_{i_j}$ such that
$|f(r(J_{i_j})) - f(\ell(J_{i_j}))| > 2^{-1} \lambda S_{i_j}$
and for all $z\in J_{i_j}$,  $f(z) \ge \min\{f(\ell (J_{i_j})), f(r (J_{i_j}))\} -
2^{-1}|f(r (J_{i_j})) - f(\ell (J_{i_j}))|$,
satisfies $\sum_j \lambda S_{i_j} > \frac2{3} \sum_{i=1}^R\lambda S_{i}$.
Call $J_{i_j}$ order preserving if
$f(\ell(J_{i_j}))< f(r(J_{i_j}))$
and order reversing if $f(\ell(J_{i_j}))> f(r(J_{i_j}))$.
We consider two subcases. In the first,  the
subcollection of order reversing intervals, which we rename as
$\{A_{1}, \dots , A_{Q}\}$
is large: $\sum_{j=1}^Q \lambda Mf(A_{j}) >
\frac1{3} \sum_{i=1}^R\lambda S_{i}$. We add the points
$\{\ell( A_{1}),r(A_{1}) \dots ,\ell (A_{Q}), r(A_{Q})\}$
to the partition $P$
and note that with this refinement the variation increases  by more than
$\frac{1}{3}\sum^R_1 \lambda S_i$. In the second subcase,
  the
subcollection of order preserving intervals,
which again we denote by
$\{A_{1}, \dots , A_{Q}\}$,
satisfies $\sum_{j=1}^Q \lambda Mf(A_{j}) >
\frac1{3} \sum_{i=1}^R\lambda S_{i}$.
Since for every $w\in A_{1}$, $f(w)
\ge f(\ell (A_{1})) - (f(r (A_{1})) - f(\ell (A_{1})))/2$,
there exists a  $c_1\in [r(A_{1}), 2r(A_{1}) - \ell (A_{1})]$
such that $f(c_1) < (f(r(A_{1})) + f(\ell (A_{1})))/2$. Otherwise, we
would have that
$$f(\ell (A_{1})) \le
\frac{1}{2(r (A_{1})-\ell (A_{1}))}
\int_{\ell (A_{1})}^{2 r (A_{1})-\ell (A_{1})} f,
$$
contradicting the assumption
that $N\subset E^c$.
Recall that
$2r(A_{1}) - \ell(A_{1}) < x_{i_1 + 1}$, so $c_1\in I_1$.
Add $\ell (A_{1}), r (A_{1})$, and $c_1$ to $P$,
together with  the endpoints of all
 intervals
$\{A_{2},\dots ,A_{{n(1)}}\}\subset \{A_{1}, \dots , A_{Q}\}$
 contained in $[r(A_1), c_1]$ (if there is any). Then
\begin{equation}
    \label{mingle}
V(f,P) + \frac12 \lambda Mf( A_{1})  +
\sum_2^{n(1)}\lambda Mf( A_{s}) < V(f, P^\prime).
\end{equation}
Go  to the next $A_{q}$  not
 already considered and repeat the process, relabeling the points in
increasing order if
needed. It may happen that  $\ell (A_{q}) < c_1 < r(A_{q})$, so
 $\ell (A_{q})$ is added
 to the left of the point $c_1$, already in the partition. But this does not
 harm any estimate, since the intervals $Mf(A_1)$ and $Mf(A_q)$ are
disjoint: If $Mf(A_q)$
lies below $Mf(A_1)$, by considering the points
$\ell (A_{1}), r (A_{1})$ and $\ell (A_{q})$ it is easily seen that the summand
$\frac12 \lambda Mf( A_{1})$ in (\ref{mingle}) can be replaced by $\lambda Mf(
A_{1})$, while if $Mf(A_q)$ lies above $Mf(A_1)$, then it is more advantageous,
from the viewpoint of guaranteeing the increase in the variation, to have  $\ell
(A_{q}) < c_1$ instead of $c_1<\ell (A_{q})$. After a finite number of steps the
list $\{A_{1}, \dots , A_{Q}\}$ is exhausted, and we get
$$
V(f,P) + \frac16
\sum_1^{R}\lambda S_{j} < V(f,P) + \frac12
\sum_1^{Q}\lambda Mf( A_{j}) < V(f, P^\prime).
$$

So regardless of whether we are in case  1 or case 2,  by (\ref{size}) we always obtain a new partition $P^\prime$ of with
$$
V(f, P) +
\frac{1}{7}\lambda^* Mf (N\cap I_1) <
V(f, P^\prime).
$$

Since all the points in $P^\prime\setminus P$ have
been chosen within $I_1$, we can repeat the argument
with every other interval $J$ determined by the first partition
$\{x_1,\dots ,x_L\}$, for which $\lambda^* Mf (N\cap J) > 0$.
In this way, a  refinement
$\{y_1,\dots ,y_K\}$
of $\{x_1,\dots ,x_L\}$ is produced, such that
$$
\sum_2^L |f(x_i) - f(x_{i-1})| + c < \sum_2^K |f(y_i) - f(y_{i-1})|.$$
\end{proof}

\section{examples.}

This section presents several examples
in order to illustrate some of the issues involved
and why different assumptions in the preceding
results are needed.

\begin{example} {\em There exists an upper semicontinuous function $f$
(with unbounded
variation) such that
$Mf$ is not continuous.}

{\em Proof:}  Let $f$ be the
characteristic function of
the closed set $\{0\}\cup \bigcup_{n=0}^\infty [3/2^{n+2}, 1/2^n]$.
Then   $Mf(0)\le 1/2$ while $\limsup_{x\to 0} Mf(x) =1 $,
so $Mf$ is discontinuous
at $0$.
\qed

This example  shows that
the hypothesis $f(w) \le Mf(w)$ in Lemma \ref{Mcont} is necessary.
We also mention that the canonical representative of $f$ is not upper
semicontinuous even though $f$ is.
\end{example}

In the next example we will follow the convention of identifying
a set with its characteristic function, thereby  using the same symbol to denote both.

\begin{example}\label{cantor} {\em There exists a bounded,
upper semicontinuous function $f$  with compact support (and unbounded
variation) such that
$Mf$ is not differentiable on a set of positive measure.
In particular, $Mf$ is not
of bounded variation.}

{\em Proof:}  We shall show that there exists a
fat Cantor set $C$ such that $MC$ is not differentiable
at any point of $C\cap \{MC
=1\}$. By a Cantor set $C$ we mean an extremely disconnected
compact set such that
all its points are points of accumulation.
Since $C$ is closed, its characteristic
function $C$ is upper semicontinuous, and obviously of unbounded variation. Note
that $MC(x) =1$ for almost every $x\in C$,  so if $MC$ is differentiable at any
such $x$, we must have $DMC(x) = 0$. If fact, $C$ will be chosen so that on $C\cap
\{MC =1\}$ the difference quotients diverge in modulus to infinity.

By fat we mean of positive measure. We shall construct $C\subset [0,1]$ so that
$\lambda C > 2/3$.  The main difference with the
usual Cantor set is that instead of removing the ``central part" of every
interval at each stage, we
remove several parts. Let $F_0 = [0,1]$ and let
$F_n$ be the finite union of closed subintervals of $[0,1]$
obtained at step $n$ of the construction, to be described below. As usual
 $C:= \cap _n F_n$. Obviously
$MC \le MF_n$; the function  $MF_n$ is the one we will
actually estimate. At stage $n$ we remove
the proportion $2^{-2n}$ of mass from the preceding set, i.e., $\lambda F_n =
(1-2^{-2n}) \lambda F_{n-1}$, so $\lambda F_{1} = 3/4$, $\lambda F_{2}= 45/64$,
 et cetera.
Then  $\lambda C = \lim_n \lambda F_{n} =1 - \lim_n \lambda F_{n}^c
> 1 - \sum_{n=1}^\infty  2^{-2n} = 2/3$. Next we ensure that for
each $n$ ``mass" and ``gaps" are sufficiently mixed. Let $I_{n-1}$
be a component of $F_{n-1}$. Subdivide $I_{n-1}$ using the $2^{2n} +
1$ equally spaced points $\ell (I_{n-1}) = x_1, x_2, \dots
,x_{2^{2n}+1}= r(I_{n-1})$, and then remove the $2^{2n} +1$ open
intervals $O(n, x_{i})$ of length $2^{-4n}\lambda I_{n-1}$, centered
at each $x_i$, noting that the first and last intervals deleted lead
only to a decrease in mass of $2^{-4n-1} \lambda I_{n-1}$ each.   Do
the same with the other components of $F_{n-1}$ to obtain $F_n$.
Then all subintervals left have the same length. Clearly, the
largest average at the points $x_i$ is obtained by considering
intervals as large as possible but without intersecting any other
deleted interval $O$, so
$$
MF_n (x_i) \le \frac{2^{-2n}
\lambda I_{n-1} - 2^{-4n} \lambda I_{n-1}}{(2^{-2n} - 2^{-4n -1}) \lambda
I_{n-1}} < 1-2^{-2n-1}.
$$
Fix $z\in C\cap \{MC =1\}$. For each $n$, let
$I_{n,z}$ be the component of $F_n$ that contains $z$,
and let $w_n$ be midpoint of the nearest  interval $O(n)$ deleted at step $n$
(if there are two such midpoints, pick any).
Then $|z - w_n|< 2^{-n(n+1)}$, since the number of components of the
 set $F_n$ is $\prod_{i=1}^n 2^{2i}=  2^{n(n+1)}$. Therefore
$$
\limsup_{w\to z} \left|\frac{MC(z) - MC(w)}{z-w}\right| \ge
\limsup_{n\to \infty} \frac{1 - MC(w_n)}{|z-w_n|}
\ge \lim_{n\to \infty} \frac{2^{-2n-1}}{2^{-n(n+1)}}
= \infty.
$$ \qed
\end{example}

\begin{remark} At present it is not clear to us how $Mf$ behaves
in higher dimensions.
This was asked in \cite{HaOn} for $W^{1,1}(\Bbb R^d)$ (Question  1, p. 169). Note that when $d>1$,
a function $f$ of bounded variation need no longer be bounded, it may not have
an upper semicontinuous representative, and
$Mf$ may fail to be continuous even if $f$ is bounded, of bounded variation, and
upper semicontinuous, as the next example shows. Furthermore the equivalence
between the pointwise variation
$V(g)$ of a function g and the size $|Dg| (\Bbb R^d)$ of its
distributional derivative no longer holds; in fact $V(g)$ is essentially a one dimensional
object, and there is no corresponding notion for $d>1$. All of this means that even if the results in dimension one
continue to hold when $d >1$,  no straightforward
extension of the arguments presented here is possible.
\end{remark}

\begin{example} {\em There exists a bounded upper semicontinuous function $f\in
BV(\Bbb R^2)$
  such that
$Mf$ is not continuous.}

{\em Proof:}  Let $f$ be the
characteristic function of
the closed triangle with vertices at $(0,0), (1,2)$, and $(2,1)$.
Then it is easy to check that
the noncentered maximal function associated to any $\ell_p$ ball
has a point of discontinuity at the origin. The same happens if
we consider the (noncentered) strong maximal operator (where averages
are taken over rectangles with sides parallel to the coordinate axes).
\qed
\end{example}

\begin{remark} Standard applications of the maximal function in
the context of $L^p$
spaces, for $p>1$, use the fact that $Mf$ dominates $|f|$ pointwise, and hence in
norm, but the norm of $Mf$ is not much larger than that of $f$. While
the latter fact is still true in Sobolev spaces by Kinnunen's theorem, $Mf$ may
fail to control $f$ in norm, as the next example shows.
This points out to the fact that applications of
$Mf$ in the theory of Sobolev spaces will tend to differ
from the usual ones in $L^p$.
One such application, explored below, consists in trying to replace $Df$ by $DMf$
in inequalities involving a function and its derivatives. Here having a smaller
norm may in fact be advantageous.
\end{remark}

\begin{example} {\em For $1\le p\le \infty$ there exists an $f\in W^{1,p}((0,1))$
  such that
$\|Mf\|_{ W^{1,p}((0,1))}<\|f\|_{ W^{1,p}((0,1))}$.}

{\em Proof:}  Let $N>>1$, set $f(x) = 1$ for $x\in (0, 2^{-1} - N^{-1})\cup (2^{-1} + N^{-1},1)$,
$f( 2^{-1})=0$, and extend $f$ piecewise linearly
to a continuous function on $(0,1)$.
Then $f$ works as advertised, since $Mf$ is ``close to being constant" (more
precisely, $Mf$ is constant save on the middle interval of length $2 N^{-1}$, where
it is Lipschitz:
 Lip$(Mf) \le \frac{1}{2^{-1} - N^{-1}}$ by Lemma
\ref{Lip})  and
$\|Mf - f\|_p$ is close to zero for
$1\le p <\infty$ (making
$N$ depend on $p$), while $\|Mf - f\|_\infty < 1$.
\qed
\end{example}

The preceding example can easily be modified
to obtain the same result in $ W^{1,p}(\Bbb R)$ for $p>1$. Also, by fixing $p$ and
letting $N$ go to infinity, one obtains a
sequence $\{f_N\}$ with  $\|Mf_N\|_{ W^{1,p}}\le c <\infty$  and
$\|f_N\|_{
W^{1,p}}\uparrow\infty$. So there is no uniform domination of $f_N$ by
$Mf_N$, even up to a constant.

\section{Applications.}
While the maximal function is a tool of every day use within the
real variable methods in harmonic analysis, its importance in the
theory of differential equations and Sobolev spaces has been considerably smaller.
This may start to change as the regularity properties of the maximal function
are being uncovered. Here we use our main result to prove
inequalities involving derivatives under less regularity, a novel kind of
 application. It is convenient in the context of Landau's inequality to adopt the
convention
$\infty
\cdot 0 = \infty$, (otherwise if $u$ is unbounded and
 $u'$ is constant, the right hand side of the inequality below is undefined).
For the real line, the sharp Landau inequality states that
given an absolutely continuous function $u^\prime$,
\begin{equation*}
  \|u'\|_\infty^2\le 2\|u\|_\infty\|u''\|_\infty.
\end{equation*}

Nowadays Landau's inequality (later generalized by Kolmogorov by considering higher order derivatives) can be regarded, save for the issue of best constants, as a special
case of the Gagliardo-Nirenberg
inequalities.
 Pointwise estimates in Landau's (and Kolmogorov's) inequality
involving the maximal function or some variant of it are known, cf. \cite{Ka},
\cite{MaSh1}, and \cite{MaSh2}. Here we present a norm inequality, involving the
derivative of the maximal function rather than the maximal function of the derivative.
As usual, $f^\prime$ denotes the  pointwise derivative of a function $f$, while $f^+:=\max\{f,0\}$ and $f^-:=
\max\{-f,0\}$ stand for its positive and negative part respectively.
We shall use $f^\prime$ and $Df$ indistinctly when $f$ is absolutely continuous.
Note that $Mf$ may be constant even if $f$ is bounded, nonnegative, and not constant.

\begin{theorem} \label{Land} Let $I$ be an interval with infinite
length, and let
$u:I\to \mathbb{R}$ be an absolutely continuous function such that
$V(u')< \infty$. Then
  \begin{equation}
    \label{landau}
    \|u'\|_\infty^2 \le 48 \|u\|_\infty \left( \|DM(u'^+)\|_\infty +
    \|DM(u'^-)\|_\infty\right).
\end{equation}
If $I =\mathbb{R}$, then
  \begin{equation}
    \label{landaur}
    \|u'\|_\infty^2 \le 24 \|u\|_\infty \left( \|DM(u'^+)\|_\infty +
    \|DM(u'^-)\|_\infty\right).
\end{equation}
\end{theorem}
\begin{proof}
 Suppose $I$ has infinite length, $\|u\|_\infty<\infty$, and
  $\|u'\|_\infty>0$.
  We claim that for every $t>0$,   \begin{equation}
    \label{landauprin}
    \|u'\|_\infty \le \max\left\{\frac{8}{t}\|u\|_\infty,
    6t\left(\|DM(u'^+)\|_\infty+\|DM(u'^-)\|_\infty\right)\right\}.
  \end{equation}
It follows, by letting $t\to\infty$, that $\|DM(u'^+)\|_\infty+\|DM(u'^-)\|_\infty>0$.
 Setting
  $$
t=\left(\frac{4\|u\|_\infty}{3(\|DM(u'^+)\|_\infty+\|DM(u'^-)\|_\infty)}\right)^{\frac12},$$
we obtain  (\ref{landau}).

To prove the claim we distinguish two cases, depending on which term of the right hand side
controls the left hand side. We may assume that $I$ is closed (otherwise we just
extend $u$ to the closure of $I$ using uniform continuity). Fix $t>0$ and
  $\alpha\in (5/6,1)$. Select $x_0 \in I$ such that
  $\max\{M(u'^+)(x_0),M(u'^-)(x_0)\} \ge \alpha \|u'\|_\infty$.
  Without loss of generality we may assume that
$ M(u'^+)(x_0)\ge \alpha \|u'\|_\infty.$

  {\em Case 1.} Suppose there exists a $y \in [x_0-t,x_0+t]\cap I$ such that
  $M(u'^+)(y) \le \frac{5}{6} \|u'\|_\infty$. By Theorem \ref{main}, $M(u'^+)$
  is absolutely continuous, so
$$
    |M(u'^+)(x_0)-M(u'^+)(y)| = \left|\int_{I(x_0,y)} DM(u'^+)\right|\le
    \|DM(u'^+)\|_\infty|x_0-y|.
$$
Hence  we have
\begin{equation}\label{landau1}
\|DM(u'^+)\|_\infty \ge \left|
\frac{M(u'^+)(x_0)-M(u'^+)(y)}{x_0-y}\right|\ge
\frac{(\alpha-5/6)\|u'\|_\infty}{t}.
\end{equation}
 {\em  Case 2.} For all $y \in [x_0-t,x_0+t]\cap I$ we have
  $M(u'^+)(y)>\frac{5}{6}\|u'\|_\infty$, so there exist  $a,b\in\Bbb R$
with $a < y < b$ such that
$$
    \frac{1}{\lambda (I\cap (a,b))}\int_{I\cap (a,b)} u'^+ \ge \frac{5}{6} \|u'\|_\infty.
$$
Write $I_y:= I\cap (a,b)$ (of course, $a$ and $b$ depend on $y$).
 Then
 $$
    \lambda(I_y\cap\{ u'^+=0\})\le \frac{1}{6}\lambda(I_y).
$$
  Now $\{I_y: y\in [x_0-t,x_0+t]\cap I\}$ is (in the subspace topology of $I$)
an open
cover  of the compact interval
  $[x_0-t,x_0+t]\cap I$, so the latter set has a finite subcover  $\{I_1, I_2,\ldots,I_N\}$. By further refining the collection, if needed, we may assume
that for every $x \in \cup_1^N I_i$,
$$
    1\le \sum_{i=1}^N \chi_{I_i}(x) \le 2
$$
(if a point belongs to three intervals, at least one of them is contained in the
union of the other two, so discard it).
 Then
  \begin{multline}\label{landau6}
    2\|u\|_\infty \ge  \int_{\cup_1^N I_i}
    u' = \int_{\cup_1^N I_i} u'^+-\int_{\cup_1^N I_i}u'^-
    \ge
    \frac{1}{2}\sum_{i=1}^N\int_{I_i}u'^+
    -\sum_{i=1}^N\int_{I_i}u'^-\\
\ge \sum_{i=1}^N\left(\frac{5}{12}\lambda\left(I_i\right)-
\lambda\left(I_i\cap \{ u'^+=0\}\right)\right)\|u'\|_\infty
    \ge \frac{1}{4}\|u'\|_\infty
    \sum_{i=1}^N\lambda(I_i)\\
\ge
     \frac{1}{4}\|u'\|_\infty
    \lambda(\cup_1^N I_i)\ge \|u'\|_\infty
    \frac{t}{4}.
  \end{multline}
 To obtain (\ref{landauprin}), put together  (\ref{landau1}) and (\ref{landau6}); then let $\alpha\to 1$.

If $I=\Bbb R$, the same argument yields $ \|u'\|_\infty
    \frac{t}{2}$,  instead of $\|u'\|_\infty
    \frac{t}{4}$, as the rightmost term in (\ref{landau6}).
 It is easy to check that this in turn gives (\ref{landaur}).
\end{proof}

 For completeness, we state the corresponding result when $I$ is bounded.
 In this case (\ref{landauprin}) is replaced by
  \begin{equation*}
    \|u'\|_\infty \le \max\left\{\frac{8}{\min\{\lambda(I),t\}}\|u\|_\infty,
    6t\left(\|DM(u'^+)\|_\infty+\|DM(u'^-)\|_\infty\right)\right\}
  \end{equation*}
for every $t>0$. Hence we obtain the following

 \begin{theorem}  Let $I$ be a bounded interval, and let
$u:I\to \mathbb{R}$ be an absolutely continuous function such that
$V(u')< \infty$. If
$$
\lambda(I)\ge
    \sqrt{\frac{4\|u\|_\infty}{3(\|DM(u'^+)\|_\infty+\|DM(u'^-)\|_\infty)}},
$$
we have
  \begin{equation*}
    \|u'\|_\infty^2 \le 48 \|u\|_\infty \left( \|DM(u'^+)\|_\infty +
    \|DM(u'^-)\|_\infty\right),
  \end{equation*}
 while if
$$
\lambda(I) <
    \sqrt{\frac{4\|u\|_\infty}{3(\|DM(u'^+)\|_\infty+\|DM(u'^-)\|_\infty)}},
$$
    we get the estimate
    \begin{equation*}
      \|u'\|_\infty \le \frac{8}{\lambda(I)}\|u\|_\infty.
    \end{equation*}
\end{theorem}

\begin{example}

Working with $DMf$ rather than with $MDf$ may lead to much better
bounds, as it happens in the following example. Let $f: \Bbb
R\to\Bbb R$ be the characteristic function of $[0,1]$. Then
$f^\prime = 0$ a.e. and  $Df = \delta_0 -\delta_{1}$, so $|Df|(\Bbb
R) = 2$. It is easy to check that $Mf (x) = x^{-1}$ if $x\ge 1$, $Mf
(x) = 1$ if $0\le x\le 1$, and $Mf (x) = (1-x)^{-1}$ if $x\le 0$, so
$Mf$ is not just Lipschitz, but even better: $DMf\in BV(\Bbb R)$. As
a side remark, we mention that for this $f$ we have $|Df|(\Bbb R) =
\|DM(f)\|_{L^1(\mathbb{R})} $, so the inequality in Theorem
\ref{main} is sharp for $\mathbb{R}$, and in fact, for every other
interval $I$: While equality is only  achieved on the real line, the
constant $1$ can never be improved, as can be seen by considering
the characteristic function of a subinterval $J\subset I$, and then
letting $J$ shrink to a point in the interior of $I$.

Now, taking $F(x):= \int_{-\infty}^y f(x) dx$, a classical Landau
inequality $(p=\infty)$ for $F$, $F^\prime$ and $F^{\prime\prime}$
would fail, since $\|F^{\prime\prime}\|_\infty =0$. Replacing
$F^{\prime\prime}$ by $DF^\prime$ does not help either, as
$\|DF^\prime\|_\infty$ makes no sense, and a natural definition
using regularizations would lead to $\|DF^\prime\|_\infty=\infty$
(in which case the inequality would be true but not useful).

Trying to extend to the setting of functions of bounded variation, the pointwise,  maximal function versions of
Landau's inequality due to Agnieszka Ka\l amajska, and
independently to  Vladimir Maz\'{}ya and  Tatyana
Shaposhnikova,  would face a similar
 difficulty:  If  the distributional derivative has a singular part,
i.e., if the function is of bounded variation but not absolutely continuous, then
its  maximal function will blow up somewhere. In the example we are considering, it
is easy to see that $M |DF^\prime|(x)
\ge \max\{|x|^{-1}, |x-1|^{-1}\}$.

Even ignoring regularity issues and replacing
$f$ with mollified versions of it, or piecewise linear continuous variants, the bounds obtained by considering $DMF^\prime$
instead of $M |DF^\prime|$ are distinctly better: Let $f_n=1 $ on $ [0,1]$, $f_n=0$ on
$(-\infty, -n^{-1})\cup (1+ n^{-1},\infty)$, and extend $f_n$ linearly in each of the
two remaining intervals, so that $f_n$ is continuous. Now $\|f_n\|_\infty =1$ and
 Theorem \ref{Land} does indeed give a bound uniform in $n$. However, the bounds obtained via
 the classical Landau inequality deteriorate as $n\to \infty$, and the same happens
 on small neighborhoods of $0$ and $1$ with pointwise inequalities using $Mf_n^\prime$.

While the function $F$ above is bounded, it is not integrable.
To obtain a similar example in
$L^1\cap L^\infty$,
 let $g(x):= f(x)-f(x-1)$ and let $G(y):=\int_{-\infty}^y g(x) dx$.
For a continuous
 example, with singular (and continuous, as a measure)
distributional derivative, let $f$ be the standard Cantor function
on $[0,1]$, defined using the Cantor ``middle third" set $C$.  We
extend it to
$\Bbb R$, first by setting $f(x) = 0$ if
$x<0$ and $f(x) = 1$ if $1\le x\le 3/2$.
Then  reflect about the axis $x=3/2$.
Now let $g(x) := f(x) - f(x-3)$. Then $G(y):=\int_{-\infty}^y g(x)
dx$ belongs to $L^1\cap L^\infty$, and, applying the notation
of Lemma \ref{Lip} to $M(g_+)$, it is clear that $\Bbb R\subset E_1$,
so $\|DMg_+\|_\infty \le 1$, and the same happens with $M(g_-)$.
Thus,
 $\|DM(g_+)\|_\infty + \|DM(g_-)\|_\infty \le 2$ and once more we
obtain a finite bound on the right hand side of inequality
(\ref{landaur}).
\end{example}

\begin{remark} It is natural to enquire whether for some constant $c$ the
simpler inequality
  \begin{equation}
    \label{landaufalse}
    \|u'\|_\infty\le c\|u\|_\infty\|DM(u')\|_\infty
  \end{equation}
holds.
  Fix $c > 0$ and select $N\in \mathbb{N}$ such that
$1> c/N$. For $k=0,\ldots,N-1$, set $u'=1$ on the intervals
  $(k/N,(2k+1)/(2N)]$, $u'=-1$ on the intervals
  $((2k+1)/(2N),(k+1)/N]$,  and $u'=0$ off $(0,1]$.
  Then $\|u'\|_\infty=1$,
   $\|u\|_\infty=1/2N$, and $|u'|=\chi_{(0,1]}$, so
  $\|DM(u')\|_\infty=1$. Thus, (\ref{landaufalse}) fails.
  \end{remark}

\begin{remark}\label{class} As indicated in the introduction, inequality
(\ref{landau}) implies Landau's, though not with the sharp constant
(in fact, the constant is not even close; the point of course
is that (\ref{landau}) can yield nontrivial bounds for some non-Lipschitz,
even discontinuous
$u'$).

Kinnunen showed that for $f\in W^{1,\infty}({\Bbb R^d})$,
$\|DMf\|_\infty \le \|Df\|_\infty$ (cf. \cite{Ki}, pages 120 and 121). The same
holds for Lipschitz functions on $I\subset \Bbb R$, as
we note in the next theorem.  While the argument is basically
the same, formally this theorem does not follow from Kinnunen's result
since we are considering also the local case $I\ne \Bbb R$.
So we include the proof. Now we have that
if $u'$ is absolutely continuous on an unbounded interval $I$
and $\|Du'\|_\infty < \infty$, then
$\|DM(u'^+)\|_\infty  + \|DM(u'^-)\|_\infty \le
\|D u'^+\|_\infty  + \|D u'^-\|_\infty \le
2\|D u'\|_\infty$, so by (\ref{landau}),
$\|u'\|_\infty\le 96 \|u\|_\infty\|Du'\|_\infty$.
  \end{remark}

\begin{theorem}  If $f:I\to \Bbb R$ is Lipschitz, then
so is
 $Mf$, and $\operatorname{Lip}(Mf) \le \operatorname{Lip}(f)$, or
equivalently
$\|DMf\|_\infty \le \|Df\|_\infty$.
The same holds for $M_Rf$.
\end{theorem}
\begin{proof} Select  $x, y\in I$, and let $f\ge 0$. Suppose $Mf(x) > Mf(y)$ and $x<y$ (in the case $x >y$ the argument is entirely analogous). If $Mf(x) =f(x)$, then
$|Mf(x) - Mf(y)|\le |f(x) - f(y)|\le \operatorname{Lip} (f) |x - y|$. Otherwise,
$Mf(x) > f(x)$, so
$$
Mf(x) = \sup_{\{[a,b]\subset I: a<b, a \le x\le b < y\}}
\frac{1}{b-a}\int_a^b f.$$ Now
$$
\frac{Mf(x) - Mf(y)}{y-x}\le
\frac{\sup_{\{[a,b]\subset I: a<b, a \le x\le b < y\}}\left(
\frac{1}{b-a}\int_a^b f -
 \frac{1}{b-a}\int_{a+ y-b}^y f\right)}{y-x}
  $$
  $$
= \sup_{\{[a,b]\subset I: a<b, a \le x\le b < y\}}
\frac{1}{(y-x)(b-a)}\int_a^b
  \left(f(t) - f(t+y-b)\right) dt
  $$
  $$
  \le \sup_{\{[a,b]\subset I: a<b, a \le x\le b < y\}}
\frac{1}{(y-x)(b-a)}\int_a^b (y-b)
  \operatorname{Lip} (f )dt
\le \operatorname{Lip} (f).
$$
 \end{proof}

In fact, the constant  $1$  given in the preceding theorem
   is not sharp.
After this paper was completed, in  joint work with Leonardo Colzani
the  authors have found the best constants:  $\operatorname{Lip}
(Mf) \le 2^{-1}\operatorname{Lip} ( f)$ for arbitrary intervals,
while on $\mathbb{R}$, $\operatorname{Lip} (Mf)  \le (\sqrt{2} - 1)
\operatorname{Lip} (f)$.
 So when deriving the classical Landau inequality from the generalization presented
 here, the constant $96$ at the end of Remark \ref{class} can be lowered
 to $48$.

To finish, we present, again under less regularity,
a trivial variant of the classical
 Poincar\'e-Wirtinger inequality, which
 states that if $f:[a,b]\to \Bbb R$ is an
absolutely continuous  function with $f(a)=f(b)=0$, then
  \begin{equation*}
    \int_a^b f(x)^2 dx \le c \int_a^b f'(x)^2dx,
  \end{equation*}
  where $c$ depends only on $b-a$.
  Using the local maximal operator $M_R$, we prove a variant of the above
  inequality, for functions of bounded variation with support at positive distance
from the boundary.
  \begin{theorem}
    Let  $f:[a,b]\to \Bbb R$ be such that  $V(f)<\infty$ and
$\operatorname{supp} f \subset
    [a+R,b-R]$ for some $R>0$. Then
    \begin{equation*}\label{MRpowi}
      \int_a^b f(x)^2 dx \le c \int_a^b DM_Rf(x)^2 dx.
      \end{equation*}
  \end{theorem}
\begin{proof} By Theorem \ref{main}, $M_Rf$ is absolutely continuous since $V(f)
<\infty$, and by hypothesis, $M_R f(a)= M_R f(b) = 0$, so we can apply the classical
Poincar\'e-Wirtinger inequality to $M_Rf$, obtaining
\begin{equation*}\label{powi}
    \int_a^b f(x)^2 dx \le  \int_a^b M_Rf(x)^2 dx \le c \int_a^b (DM_Rf(x))^2dx,
  \end{equation*}
\end{proof}

\end{document}